\documentclass[pdflatex,sn-mathphys-num]{sn-jnl}
\usepackage{listings}

\usepackage{amsmath}
\usepackage{amsfonts}
\usepackage{bm}
\usepackage{enumerate}
\usepackage{graphicx}
\usepackage{float} 
\usepackage{subfigure} 
\usepackage{caption} 
\usepackage[section]{placeins}
\usepackage{mathrsfs}
\usepackage{subfigure}
\usepackage{url}
\usepackage{verbatim}
\usepackage{amssymb,graphicx}
\usepackage{tabularx}
\usepackage{makecell}

\theoremstyle{thmstyleone}%
\theoremstyle{thmstyletwo}%

\theoremstyle{thmstylethree}%
\newtheorem{definition}{Definition}%

\begin{document}

\title[Semi-implicit-explicit Runge-Kutta method for nonlinear differential equations]{Semi-implicit-explicit Runge-Kutta method for nonlinear differential equations}

 \author[1]{\fnm{Lingyun} \sur{Ding}}\email{dingly@g.ucla.edu}

\affil[1]{\orgdiv{Department of Mathematics}, \orgname{University of California, Los Angeles}, \orgaddress{\city{Los Angeles}, \postcode{90024}, \state{CA}, \country{USA}}}

\abstract{
A semi-implicit-explicit (semi-IMEX) Runge-Kutta (RK) method is proposed for the numerical integration of ordinary differential equations (ODEs) of the form $\mathbf{u}' = \mathbf{f}(t,\mathbf{u}) + G(t,\mathbf{u}) \mathbf{u}$, where $\mathbf{f}$ is a non-stiff term and $G\mathbf{u}$ represents the stiff terms. Such systems frequently arise from spatial discretizations of time-dependent nonlinear partial differential equations (PDEs). For instance, $G$ could involve higher-order derivative terms with nonlinear coefficients. Traditional IMEX-RK methods, which treat $\mathbf{f}$ explicitly and $G\mathbf{u}$ implicitly, require solving nonlinear systems at each time step when $G$ depends on $\mathbf{u}$, leading to increased computational cost and complexity. In contrast, the proposed semi-IMEX scheme treats $G$ explicitly while keeping $\mathbf{u}$ implicit, reducing the problem to solving only linear systems. This approach eliminates the need to compute Jacobians while preserving the stability advantages of implicit methods. A family of semi-IMEX RK schemes with varying orders of accuracy is introduced. Numerical simulations for various nonlinear equations, including nonlinear diffusion models, the Navier-Stokes equations, and the Cahn-Hilliard equation, confirm the expected convergence rates and demonstrate that the proposed method allows for larger time step sizes without triggering stability issues.
}



\maketitle
\section{Introduction}
This paper addresses the numerical integration of ordinal differential equations (ODEs) of the form
\begin{equation}\label{eq:vector evolution equation}
\mathbf{u}' = \mathbf{f}(t,\mathbf{u}) + G(t,\mathbf{u}) \mathbf{u},
\end{equation}
where $\mathbf{u}(t) \in \mathbb{R}^N$ is an unknown vector function, $\mathbf{f}: \mathbb{R}^N \times \mathbb{R} \to \mathbb{R}^N$ is a given vector field, and $G:\mathbb{R} \times \mathbb{R}^{N}\to \mathbb{R}^{N \times N}$ is a matrix. $\mathbf{f}$ and $G$ are assumed to be Lipschitz continuous with respect to $\mathbf{u}$. Systems of the form \eqref{eq:vector evolution equation} often arise from spatial discretizations of time-dependent partial differential equations (PDEs), such as Cahn-Hilliard equation~\cite{cahn1958free,hester2023fluid}, Navier-Stokes equation \cite{sanderse2012accuracy,froehle2014high},  Keller-Segel equation~\cite{zayernouri2016fractional}, thin-film equations~\cite{kondic2003instabilities,mavromoustaki2018surface}. In this context, the vector $\mathbf{u}$ represents the discretized unknown field, while $\mathbf{f}$ corresponds to a forcing term or the non-stiff part of the equation, and $G\mathbf{u}$ accounts for the stiff terms. For instance, $G$ could involve higher-order derivative terms with nonlinear coefficients.

Explicit methods are computationally efficient and easy to implement. However, in the presence of a stiff term, they require severely restricted time steps, often several orders of magnitude smaller than those used in non-stiff cases \cite{butcher1996history,quarteroni2010numerical}. Additionally, PDE systems typically involve boundary conditions that are challenging to enforce within an explicit scheme.

A popular alternative is the implicit-explicit (IMEX) Runge-Kutta (RK) method, which treats the non-stiff term $\mathbf{f}$ explicitly and the stiff term $G\mathbf{u}$ implicitly \cite{pareschi2005implicit,ascher1997implicit,boscarino2024implicit,cavaglieri2015low}. This approach improves stability by ensuring that the stability constraints are determined by the non-stiff term alone and allows for boundary conditions to be enforced. However, when $G$ depends nonlinearly on $\mathbf{u}$, each time step involves solving a nonlinear system, which can be computationally expensive and may introduce stability and sensitivity issues.

Several approaches have been proposed to mitigate this issue. One such approach is the Rosenbrock method \cite{norsett1979order,shen1996semi,albrecht1996runge}, which approximates the stiff term as the product of the Jacobian matrix and the solution vector. In this method, the Jacobian matrix is treated explicitly, while the solution vector is treated implicitly. Efforts have been made to reduce the number of Jacobian matrix evaluations \cite{steihaug1979attempt}, though occasional evaluations are still required. Specialized schemes of this type have been developed for solving the Navier-Stokes equations \cite{nikitin2006third} and multibody mechanical systems \cite{song2004semi}.

The second related approach is the semi-IMEX RK method for ODEs with stiff damping terms, as studied in \cite{chertock2015steady}. The equation considered in their work has a similar form to \eqref{eq:vector evolution equation}, but with $G$ as a diagonal matrix representing the damping term. Following a similar strategy, \cite{hu2018asymptotic} developed a scheme for the Bhatnagar-Gross-Krook kinetic equation. This class of methods typically consists of two stages: the first stage employs classical explicit schemes or IMEX schemes with appropriate modifications, while the second stage introduces a correction to enhance accuracy or enforce additional properties.

Third, some schemes are tailored to specific equations. For example, \cite{mata2011numerical} developed a second-order time-stepping scheme for the equation governing particle-laden thin film flow, only solving linear systems while avoiding the use of the Jacobian matrix. However, there is no clear procedure for generalizing this approach to other equations.

The above methods do not fully exploit the structure of equation \eqref{eq:vector evolution equation}. We propose a semi-implicit-explicit scheme that treats $G(t, \mathbf{u})$ explicitly while keeping $\mathbf{u}$ implicit. As a result, instead of solving nonlinear systems, this method only requires solving a sequence of linear systems at each time step, significantly reducing computational cost. The proposed approach is simpler and easier to implement compared to the aforementioned methods while eliminating the need to compute the Jacobian matrix and preserving the Runge-Kutta formulation.

To illustrate the approach, consider the advection-diffusion equation with a nonlinear diffusion coefficient $D(c)$ depends on the concentration:
\begin{equation}\label{eq:nonlinear advection-diffusion equation}
\partial_{t} c + \mathbf{u} (x,t) \cdot \nabla c = \nabla \cdot \left( D(c) \nabla c \right)+S (x,t),
\end{equation}
where $c$ represents the concentration field, $\mathbf{u}$ is a prescribed velocity field, and $S (x,t)$ is a source term. Such equations commonly arise in models of particle suspension \cite{griffiths2012axial,lee2025comparative,ding2025equilibrium},  stratified fluid \cite{ding2024diffusion}, porous media \cite{aronson2006porous} and  electrolyte solution \cite{gupta2019diffusion,ding2023shear}.

The second-order derivative term (diffusion term) is generally more stiff compared to the first-order derivative term (advection term). Therefore, the advection term $\mathbf{u} \cdot \nabla c$ corresponds to $\mathbf{f}$ in equation \eqref{eq:vector evolution equation}, while the diffusion operator $\nabla \cdot \left( D(c) \nabla \cdot \right)$ corresponds to $G$ in equation \eqref{eq:vector evolution equation}. The first-order IMEX-RK method employs the following time discretization (suppressing spatial discretization):
\begin{equation}
\frac{c_{n+1} - c_{n}}{h} + \mathbf{u} (x,t_{n}) \cdot \nabla c_{n} = \nabla \cdot \left( D(c_{n+1}) \nabla c_{n+1} \right)+S (x,t_{n}),
\end{equation}
which results in a nonlinear system for $c_{n+1}$. $h=t_{n+1}-t_{n}$ is the time step size. $c_{n}=c (x,t_{n})$. 

In contrast, the proposed semi-implicit-explicit method discretizes the equation as
\begin{equation}\label{eq:example ADE FB}
\frac{c_{n+1} - c_{n}}{h} + \mathbf{u} (x,t_{n}) \cdot \nabla c_{n} = \nabla \cdot \left( D(c_{n}) \nabla c_{n+1} \right)+S (x,t_{n}),
\end{equation}
where the diffusion coefficient $D(c_{n})$ is treated explicitly. As a result, the system remains linear in $c_{n+1}$, making it easier to solve compared to the nonlinear case. This scheme is first-order accurate in time. The objectives of this paper are to develop higher-order schemes, analyze their theoretical properties, and investigate their numerical implementation.

The paper is organized as follows: Section \ref{sec:Semi IMEX-RK schemes} introduces the general formulation of the semi-IMEX RK scheme and presents methods with varying accuracy orders. Section \ref{sec:numerical test} presents numerical tests for various equations (scalar equation, nonlinear diffusion equation, Cahn-Hilliard equation, Navier-Stokes equation), verifying the convergence rate and highlighting the advantages of the semi-IMEX RK scheme in specific scenarios. Finally, Section \ref{sec:conclusion} concludes the paper with several research directions.


\section{Semi-IMEX-RK schemes}
\label{sec:Semi IMEX-RK schemes}
We first review the formulation of the IMEX RK scheme, which combines an $s$-stage explicit scheme and an $s$-stage implicit scheme. It can be represented using a double tableau in the standard Butcher notation:
\begin{equation}
\begin{array}{c|cccccc}
0& 0&0&...&0&0\\
 \tilde{c}_{2}& \tilde{a}_{21}&0&...&0&0\\
\tilde{c}_{3}& \tilde{a}_{31}& \tilde{a}_{32}&...&0&0\\
...\\
\tilde{c}_{s}&\tilde{a}_{s,1}& \tilde{a}_{s,2}& ...&\tilde{a}_{s,s-1}& 0\\
\hline
&\tilde{b}_{1}& \tilde{b}_2& ...&\tilde{b}_{s-1}&\tilde{b}_{s}\\
 \end{array}
, \quad
\begin{array}{c|ccccc}
 c_{1}& a_{11}&0&...&0\\
 c_{2}& a_{21}&a_{22}&...&0\\
c_{3}& a_{31}& a_{32}&...&0\\
...\\
c_{s}&a_{s1}& a_{s2}& ...& a_{s,s}\\
\hline
&b_{1}& b_2& ...&b_s\\
 \end{array}.
\end{equation}
Here, the parameters with tildes correspond to the explicit scheme. One step of the classical IMEX-RK scheme for equation \eqref{eq:vector evolution equation} from $t_n$ to $t_{n+1}$ is given by
\begin{equation}\label{eq:IMEX-RK scheme}
  \begin{aligned}
&\left( 1- h a_{i,i}G (t_{n}+c_{i}h, K_{i}) \right)K_{i}= u_{n}+ h  \sum\limits_{j=1}^{i-1} \tilde{a}_{i,j} f (t_{n}+\tilde{c}_{j}h, K_{j}) + a_{i,j}G (t_{n}+c_{i}h, K_{j})K_{j},\\
&u_{n+1}= u_{n} + h \sum\limits_{i=1}^{s} \tilde{b}_{j}f(t_{n}+ \tilde{c}_{j}h, K_{j})+ b_{j}G(t_{n}+c_{j}h,  K_{j})K_{j},\\ 
\end{aligned}
\end{equation}
where  $ i=1,...,s$.  For $i$-th stage, we need to solve a nonlinear equation for $K_{i}$. To avoid solving nonlinear systems, we modify the above scheme as follows:
\begin{equation}\label{eq:semi-IMEX-RK scheme}
  \begin{aligned}
&\left( 1- h a_{i,i}G (t_{n}+c_{i}h, \tilde{K}_{i}) \right)K_{i}= u_{n}+ h  \sum\limits_{j=1}^{i-1} \tilde{a}_{i,j} f (t_{n}+\tilde{c}_{j}h, K_{j}) + a_{i,j}G (t_{n}+c_{i}h, K_{j})K_{j},\\
&u_{n+1}= u_{n} + h \sum\limits_{i=1}^{s} \tilde{b}_{j}f(t_{n}+ \tilde{c}_{j}h, K_{j})+ b_{j}G(t_{n}+c_{j}h,  K_{j})K_{j}+ b_{i+1}G (t_{n}+c_{i}h, \tilde{K}_{i}) K_{i}.\\ 
\end{aligned}
\end{equation}
There are two modifications. The first modification is replacing $a_{i,i}G (t_{n}+c_{i}h, K_{i})$ in $i$-th equation with $a_{i,i}G (t_{n}+c_{i}h, \tilde{K}_{i})$, where $\tilde{K}_{i}$ depends on $K_{j}$ for $1\leq j<i$. Therefore, each stage $K_{i}$ is obtained by solving a linear system. The simplest choice is to set $\tilde{K}_{i} = K_{i-1}$ for $i>1$ and $\tilde{K}_{1} = u_{n}$. Alternatively, $\tilde{K}_{i}$ can be approximated using an explicit scheme. Since this alternative approach can be reformulated as the first case by increasing the number of stages and appropriately selecting the method parameters, equation \eqref{eq:semi-IMEX-RK scheme} remains applicable (examples will be provided later). Therefore, to maintain a consistent formulation, we set $\tilde{K}_{i} = K_{i-1}$ unless otherwise specified. 

The second modification involves adding the additional term $b_{i+1}G (t_{n}+c_{i}h, \tilde{K}_{i}) K_{i}$ when computing $u_{n+1}$.  Since the matrix $G (t_{n}+c_{i}h, \tilde{K}_{i})$ is constructed at the $i$-th stage, incorporating this term allows for the reuse of the matrix while ensuring that the scheme satisfies certain desired properties.

We are looking for the methods that ideally satisfy the following desirable properties.  The first concerns boundary conditions or other constraints in differential equations, such as the no-slip boundary condition and incompressibility in the Navier-Stokes equations. These conditions are typically enforced when solving the linear or nonlinear system arising from the implicit scheme, ensuring that $K_s$ satisfies the boundary conditions and constraints. If $u_{n+1}$ is computed as a linear combination of $f(K_i)$ and $G (K_{i})K_{i}$, it may not inherently preserve these properties. Therefore, we are particularly interested in methods that satisfy $u_{n+1} = K_s$, which requires the conditions $b_{j}=a_{s,j}$, $\tilde{b}_{j}=\tilde{a}_{s,j}$ for $j=1,...,s-1$,  $b_{s+1}=a_{s,s}$, $\tilde{b}_{s}=0$. In some cases, for a given number of stages and accuracy order, such a method may not be available. In that case, we relax the condition to %
\begin{equation}\label{eq:alpha type method}
\begin{aligned}
\alpha b_{j}=a_{s,j}, \alpha \tilde{b}_{j}=\tilde{a}_{s,j},  j=1,...,s-1,\;  \alpha b_{s+1}=a_{s,s}, \tilde{b}_{s}=0,
\end{aligned}
\end{equation}where $\alpha$ is a constant. This results in $u_{n+1}=\frac{1}{\alpha}K_{s}+(1-\frac{1}{\alpha})u_{n}$. If the boundary condition is linear and time-independent, this ensures that $u_{n+1}$ satisfies the boundary condition.

The second property is stability. We primarily focus on the stability properties of the implicit part of the scheme, as they justify the combination of explicit and implicit methods. A-stability ensures that a method remains stable for all step sizes when applied to linear test problems with negative real eigenvalues. This property is particularly desirable for handling stiff equations without step size restrictions. The precise definition is as follows:
\begin{definition}[Stability function]
Applying scheme \eqref{eq:semi-IMEX-RK scheme} to the ODE with $f (t,u) =0, g (t,u)=\lambda$, resulting a relation between $u_{n+1}$ and $u_{n}$: $u_{n+1}=R(\lambda h)u_{n}$. Here $R (z)$ is the stability function. 
\end{definition}

\begin{definition}[A stable]
A method is A-stable if its stability function satisfies $|R(z)| < 1$ for all $z \in \mathbb{C}$ with $\Re(z) < 0$. The method is L-stable if it is A-stable and additionally satisfies $R(z) \to 0$ as $z \to -\infty$.
\end{definition}

Third, we seek methods that satisfy the row-sum conditions: $\sum\limits_{j=1}^{i-1}\tilde{a}_{i,j}=\tilde{c}_{i}$,  $\sum\limits_{j=1}^{i}a_{i,j}=c_{i}$ for $i=1,...,s$.  Although this condition is neither sufficient nor necessary for the order of accuracy, it is commonly used in the literature and can simplify the order conditions in some cases. Thus, we have decided to adopt it.

The order conditions can be derived using a Taylor expansion. In the following sections, we will present the conditions and schemes up to order 3, where the row-sum condition has been used to simplify the analysis.

\subsection{First order}
The conditions for a first-order method are given by
\begin{equation}\label{eq:order condition O1}
\begin{aligned}
&\sum\limits_{i=1}^{s+1}\tilde{b}_{i}= \sum\limits_{i=1}^{s}b_{i}=1.\\
\end{aligned}
\end{equation}

The scheme in equation \eqref{eq:example ADE FB} corresponds to the forward-backward Euler method, with the associated tableau shown in Table \ref{tab:forward backward Euler}. In this case, the scheme \eqref{eq:semi-IMEX-RK scheme} simplifies to
\begin{equation}\label{eq:forward-backward Euler}
\begin{aligned}
&u_{n+1}=u_{n}+h \left(  f (t_{n}, u_{n})+G (t_{n}+h, u_{n})u_{n+1} \right).
\end{aligned}
\end{equation}
The corresponding stability function is $ R (z) = \frac{1}{1 - z}$. Since $\lim\limits_{z \to -\infty} R(z) = 0$, the method is L-stable.
\begin{table}[h]
\centering
\begin{minipage}{0.4\textwidth}
\centering
\begin{tabular}{c|cc}
  0 & 0 & 0 \\
  1 & 1 & 0\\
\hline
  & 1 &0
\end{tabular}
\end{minipage}
\hspace{0.02\textwidth}
\begin{minipage}{0.4\textwidth}
\centering
\begin{tabular}{c|ccc}
  0 & 0 & 0 &\\
  1 & 0 & 1&\\
\hline
  & 0 &0&1
\end{tabular}
\end{minipage}
 \caption{Tableau for the explicit(left) and implicit(right) forward-backward Euler scheme.}
 \label{tab:forward backward Euler}
\end{table}

\subsection{Second order}
In addition to  condition \eqref{eq:order condition O1}, six conditions must be satisfied for the scheme to achieve second-order accuracy:
\begin{equation}\label{eq:order condition O2}
\begin{aligned}
  &  \sum\limits_{i=1}^{s}b_{i}c_{i}+b_{s+1}c_{s-1}=\frac{1}{2},\quad 
  \sum\limits_{i=2}^{s}b_{i}\tilde{c}_{i}+ b_{s+1}\tilde{c}_{s-1} =\frac{1}{2},\quad  \sum\limits_{i=1}^{s}b_{i}c_{i}+H (s-2)b_{s+1}c_{s}=\frac{1}{2}\\
&  \sum\limits_{i=1}^{s}\tilde{b}_{i}c_{i}= \frac{1}{2}, \quad   \sum\limits_{i=2}^{s}b_{i}\tilde{c}_{i}+b_{s+1}\tilde{c}_{s}=\frac{1}{2},\quad  \sum\limits_{i=2}^{s}\tilde{b}_{i}\tilde{c}_{i}  =\frac{1}{2}.
\end{aligned}
\end{equation}
where $H (x)=0$ when $x\leq0$ and 1 when $x>0$. For comparison, the IMEX RK method in \cite{pareschi2005implicit} requires two conditions for first-order accuracy and four additional conditions for second-order accuracy.

One example of second order method is the midpoint method with the table \ref{tab:midpoint}.
\begin{table}[h]
\centering
\begin{minipage}{0.4\textwidth}
\centering
\begin{tabular}{c|cc}
  0 & 0 & 0 \\
 $\frac{1}{2}$ & $\frac{1}{2}$ & 0\\
\hline
  & 0 &1
\end{tabular}
\end{minipage}
\hspace{0.02\textwidth}
\begin{minipage}{0.4\textwidth}
\centering
\begin{tabular}{c|ccc}
  0 & 0 & 0 \\
  $\frac{1}{2}$ & 0 &$\frac{1}{2}$\\
\hline
  & 0 &1
\end{tabular}
\end{minipage}
 \caption {Tableau for the explicit(left) and implicit(right) midpoint method, second order scheme.}
 \label{tab:midpoint}
\end{table}

With this tableau,  equation \eqref{eq:semi-IMEX-RK scheme} simplifies to
\begin{equation}
  \begin{aligned}
&   K_{2}= u_{n}+ \frac{h}{2} f \left( t_{n}, u_{n}\right) + \frac{h}{2}G \left( t_{n}+\frac{h}{2}, u_{n} \right)K_{2},\\ 
&u_{n+1}= u_{n} +h \left(  f \left( t_{n}+\frac{h}{2},K_{2} \right) + G \left( t_{n}+\frac{h}{2}, K_{2} \right)K_{2} \right) .\\
\end{aligned}
\end{equation}
The stability function is $R (z)=\frac{2+z}{2-z}$. So the method is A-stable. 

The three-stage, second-order methods satisfying condition \eqref{eq:alpha type method} admit the following tableau \ref{eq:three stage second order alpha type method}.
\begin{table}[h]
\centering
\begin{minipage}{0.4\textwidth}
\centering
\begin{tabular}{c|ccc}
    0& 0& 0 & 0\\
$\alpha$  & $\alpha$ & 0 &0 \\
$\alpha$ &$\frac{2\alpha-1}{2}$  &$\frac{1}{2}$ & 0 \\
\hline
   &$\frac{2\alpha-1}{2\alpha}$ &$ \frac{1}{2\alpha}$& 0
\end{tabular}
\end{minipage}
\hspace{0.02\textwidth}
\begin{minipage}{0.4\textwidth}
\centering
\begin{tabular}{c|cccc}
    0& 0& 0 & 0&\\
$\alpha$ &$a_{2,1}$  & $\alpha-a_{2,1}$ &0 &\\
$\alpha$ &$\frac{2 \alpha -1}{2}$   & $\frac{1-2 \alpha b_{4}}{2}$ &$\alpha b_{4}$& \\
\hline
   &$\frac{2 \alpha -1}{2 \alpha}$  &$\frac{1-2 \alpha b_{4}}{2 \alpha}$ & 0&$b_{4}$
\end{tabular}
\end{minipage}
  \caption {Tableau for the explicit(left) and implicit(right) second order scheme satisfying condition \eqref{eq:alpha type method}.}
  \label{eq:three stage second order alpha type method}
\end{table}

The stability function is
\begin{equation}\label{eq:stability function stage3 alpha type}
\begin{aligned}
&R (z)= \frac{z^2 \left(2 \alpha  b_4\left(-a_{2,1}+\alpha -1\right)+2 a_{2,1}-2 \alpha +1\right)-2 z \left(-a_{2,1}+\alpha +\alpha  b_4-1\right)+2}{2 \alpha  b_4z^2 \left(\alpha -a_{2,1}\right)-2 z \left(-a_{2,1}+\alpha +\alpha  b_4\right)+2}.
\end{aligned}
\end{equation}
We consider two special cases. First, we set $\alpha=\frac{1}{2}, b_{4}=1, a_{21}=0$ to minimize the number of nonzero elements in the tableau. The resulting tableau is table \ref{tab:trapezoid}.
\begin{table}[h]
\centering
\begin{minipage}{0.4\textwidth}
\centering
\begin{tabular}{c|ccc}
    0& 0& 0 & 0\\
$\frac{1}{2}$  &$ \frac{1}{2}$ & 0 &0 \\
$\frac{1}{2}$  &0 &$\frac{1}{2}$ & 0 \\
\hline
   &0 &1& 0
\end{tabular}
\end{minipage}
\hspace{0.02\textwidth}
\begin{minipage}{0.4\textwidth}
\centering
\begin{tabular}{c|cccc}
    0& 0& 0 & 0&\\
$\frac{1}{2}$ &0  & $\frac{1}{2}$ &0 &\\
$\frac{1}{2}$ & 0  &0 & $\frac{1}{2}$& \\
\hline
   &0  &0 & 0&1 
\end{tabular}
\end{minipage}
  \caption {Tableau for the explicit(left) and implicit(right) second order scheme.}
  \label{tab:trapezoid}
\end{table}

The corresponding scheme is
\begin{equation}
  \begin{aligned}
&   K_{2}= u_{n}+   \frac{h}{2}f \left( t_{n}, u_{n} \right) + \frac{h}{2}G \left( t_{n}+\frac{h}{2}, u_{n} \right)K_{2}, \\ 
&   K_{3}= u_{n}+ \frac{h}{2} f \left( t_{n}+\frac{h}{2}, K_{2} \right) + \frac{h}{2}G \left( t_{n}+\frac{h}{2}, K_{2} \right)K_{3},\\ 
&u_{n+1}= u_{n} +h \left(  f \left( t_{n}+\frac{h}{2}, K_{2} \right) + G \left( t_{n}+\frac{h}{2}, K_{2} \right)K_{3} \right) = 2K_{3}-u_{n}.\\
\end{aligned}
\end{equation}
The stability function is also $R (z)=\frac{2+z}{2-z}$, implying A-stable.

Second, to pursue an $L$-stable method, based on equation \eqref{eq:stability function stage3 alpha type}, we derive the necessary condition $2 \alpha  b_4\left(-a_{2,1}+\alpha -1\right)+2 a_{2,1}-2 \alpha +1=0$ along with the constraint $\alpha -a_{2,1}\neq 0$. To ensure uniform diagonal elements, we impose the condition  $\alpha -a_{2,1}=b_4$. This leads to two sets of solutions:
\begin{equation}
\begin{aligned}
&a_{21}=\frac{\pm(2 \alpha ^2-\alpha -1) \pm \sqrt{\alpha ^2+1}}{2 \alpha }, b_{4}=\frac{\mp\sqrt{\alpha ^2+1}+\alpha +1}{2 \alpha}.
\end{aligned}
\end{equation}
When $\alpha = 1$, one set of solutions falls outside the range $[0,1]$. Thus, we select $a_{21}=\frac{1}{\sqrt{2}},   b_{4}=\frac{1}{2} \left(2-\sqrt{2}\right)$. The corresponding Butcher tableau is given in Table \ref{tab:second order scheme L}.
\begin{table}[h]
\centering
\begin{minipage}{0.4\textwidth}
\centering
\begin{tabular}{c|ccc}
    0& 0& 0 & 0\\
1  & 1 & 0 &0 \\
1 &$\frac{1}{2}$  &$\frac{1}{2}$ & 0 \\
\hline
   &$\frac{1}{2}$   &$\frac{1}{2}$& 0
\end{tabular}
\end{minipage}
\hspace{0.02\textwidth}
\begin{minipage}{0.4\textwidth}
\centering
\begin{tabular}{c|cccc}
    0& 0& 0 & 0&\\
1  &$\frac{1}{\sqrt{2}}$  & $\frac{2-\sqrt{2}}{2}$ &0 &\\
1 &$\frac{1}{2}$  & $\frac{1}{\sqrt{2}}-\frac{1}{2}$& $\frac{2-\sqrt{2}}{2}$& \\
\hline
   &$\frac{1}{2}$   &$\frac{1}{\sqrt{2}}-\frac{1}{2}$ & 0&$\frac{2-\sqrt{2}}{2}$
\end{tabular}
\end{minipage}
  \caption {Tableau for the explicit(left) and implicit(right) second order L-stable scheme.}
 \label{tab:second order scheme L}
\end{table}

The second-order semi-IMEX RK method can be derived from the classical second-order IMEX RK method. With the row sum condition, $K_{i}$ can be regarded as an approximation of the solution at $t = t_{n} + c_{i}h$. We choose $\tilde{K}_{i} = u{n} + c_{i}h \big(f (t_{n}, u_{n}) + G(t_{n}, u_{n}) u_{n} \big)$. Then, $K_{i}$ in equation \eqref{eq:semi-IMEX-RK scheme} differs from $K_{i}$ in equation \eqref{eq:IMEX-RK scheme} by $\mathcal{O}(h^{3})$. Therefore, the same Butcher table can be used while preserving the second-order accuracy. In some cases, an even simpler expression for $\tilde{K}_{i}$ may be chosen. For example,  tables \ref{tab:IMEX second order} correspond to a second-order $L$-stable IMEX method from \cite{pareschi2005implicit}.
\begin{table}[h]
\centering
\begin{minipage}{0.4\textwidth}
\centering
\begin{tabular}{c|cc}
  0 & 0 & 0 \\
1 &1 & 0\\
\hline
  &  $\frac{1}{2}$ &  $\frac{1}{2}$
 \end{tabular}
\end{minipage}
\hspace{0.02\textwidth}
\begin{minipage}{0.4\textwidth}
\centering
 \begin{tabular}{c|cc}
  $\gamma$ & $\gamma$ & 0 \\
$1-\gamma$  &$1-2\gamma$& $\gamma$\\
\hline
  &$\frac{1}{2}$&$\frac{1}{2}$
 \end{tabular}
\end{minipage}
 \caption{Tableau for the explicit (left) and implicit (right) second-order IMEX Runge-Kutta scheme, where $\gamma = 1 - \frac{1}{\sqrt{2}}$. This tableau should be interpreted in the context of equation \eqref{eq:IMEX-RK scheme}.  }
\label{tab:IMEX second order} 
\end{table}
With this table, we can construct a second order $L$-stable semi-IMEX method
\begin{equation}
  \begin{aligned}
&   K_{1}= u_{n}+  \gamma h  G \left( t_{n}+\gamma h, u_{n} \right)K_{1}, \quad  \tilde{K}_{1}=u_{n}+(1-\gamma) h G \left( t_{n}+\gamma h, u_{n} \right)u_{n},\\ 
&   K_{2}= u_{n}+ h (f \left( t_{n}, K_{1} \right) + (1-2\gamma)G \left( t_{n}+\gamma h, K_{1} \right)K_{1} + \gamma  G \left( t_{n}+(1-\gamma) h, \tilde{K}_{1} \right)K_{2}),\\ 
&u_{n+1}= u_{n} +\frac{h}{2} \left( f \left( t_{n}, K_{1} \right)+  f \left( t_{n}+ h, K_{2} \right) + G \left( t_{n}+\gamma h, K_{1} \right)K_{1} + G \left( t_{n}+ (1-\gamma)h, K_{2} \right)K_{2} \right) .\\
\end{aligned}
\end{equation}
We can reformulate the above method as a three-stage scheme so that the procedure in equation \eqref{eq:semi-IMEX-RK scheme} remains applicable. Specifically, we treat $K_{1}, \tilde{K}_{1}, K_{2}$ in the above method as $K_{1}, K_{2}, K_{3}$ in the three-stage method descried by  tableau \ref{tab:second order scheme 3}.  
\begin{table}[h]
\centering
\begin{minipage}{0.4\textwidth}
\centering
  \begin{tabular}{c|ccc}
   0 & 0 & 0 &0\\
  0 & 0 & 0 &0\\
1 &  1 &0 & 0 \\
\hline
  & $\frac{1}{2}$&  0 &  $\frac{1}{2}$
 \end{tabular}
\end{minipage}
\hspace{0.02\textwidth}
\begin{minipage}{0.4\textwidth}
\centering
 \begin{tabular}{c|ccc}
   $\gamma$ & $\gamma$ & 0& 0 \\
$1-\gamma$ & $1-\gamma$ & 0& 0 \\
$1-\gamma$  &$1-2\gamma$& 0&$\gamma$\\
\hline
  & $ \frac{1}{2}$ &0& $\frac{1}{2}$
 \end{tabular}
\end{minipage}
 \caption{Tableau for the explicit(left) and implicit(right) second order scheme, where $\gamma=1-\frac{1}{\sqrt{2}}$}
  \label{tab:second order scheme 3}  
\end{table}

\subsection{Third order}
The conditions for achieving third-order accuracy are highly complex. When the row sum condition is considered, in addition to the conditions required for first- and second-order methods, there are 29 additional conditions for third order, while the IMEX method described in \cite{pareschi2005implicit} requires 14. These extra conditions arise from replacing $1 - a_{i,i} G (t_{n} + c_{i}h, K_{i}) K_{i}$ with $1 - a_{i,i} G (t_{n} + c_{i}h, K_{i-1}) K_{i}$, introducing additional couplings.

\begin{table}[h]
  \footnotesize
\centering
\begin{tabular}{c@{\hspace{0.16cm}}|c@{\hspace{0.16cm}}c@{\hspace{0.16cm}}c@{\hspace{0.16cm}}c@{\hspace{0.16cm}}}
  0 & 0& 0 & 0&0\\
  0.7775079538595848 & 0.7775079538595848& 0 &0&0 \\
0.6583867604773560 &0.3850382624054263 &0.2733484980719337 & 0&0 \\
0.6583867604773565  &0.2905474198112961 &0.1784065415104640 & 0.1894327991556034&0 \\
\hline
   &0.2486553715043413  &0.04469938464765911 &0.3828282521031255 &0.3238169917448679 
\end{tabular}
\begin{tabular}{c@{\hspace{0.16cm}}|c@{\hspace{0.16cm}}c@{\hspace{0.16cm}}c@{\hspace{0.16cm}}c@{\hspace{0.16cm}}}
  0& 0& 0 & 0&0\\
0.7775079538595848 &0.5668275181562270 &0.2106804357033578  &0 &0 \\
0.6583867604773565&0.3481097445529071 & 0.1497169356151823  &0.1605600803092672  &0 \\
 0.6583867604773565& 0.3299758037920577&0.1113697479208660 & 0.1255619659848192 &0.09147924277961349\\
\hline
   &0.2486553715043413  &0.04469938464765911 &0.3828282521031255 &0.3238169917448679 
\end{tabular}
\caption{Tableau for the explicit(top) and implicit(bottom) third order $L$-stable scheme.}
\label{tab:third order four stage}
  \end{table}

Given this complexity, we omit listing the order conditions and instead provide only the Butcher tableaus for the methods. Table \ref{tab:third order four stage} shows a fourth stage, third order method. The stability function is
 \begin{equation}
R(z)=\frac{33.95359 z^2+173.6267 z+323.1586}{-z^3+21.90616 z^2-149.5318 z+323.1586}.
\end{equation}
The method is L-stable.

Table \ref{tab:third order five stage} presents a five-stage, third-order method, stratifying $u_{n+1}=K_{s}$. Although it is formulated as a five-stage method, it requires solving only three linear systems in the update  from $t_n$ to $t_{n+1}$, making its computational cost comparable to the four-stage method in Table \ref{tab:third order four stage}.
\begin{table}[h]
  \footnotesize
\centering
\begin{tabular}{c@{\hspace{0.16cm}}|c@{\hspace{0.16cm}}c@{\hspace{0.16cm}}c@{\hspace{0.16cm}}c@{\hspace{0.16cm}}c@{\hspace{0.16cm}}}
  0 & 0& 0 & 0&0 &0\\
  0.6411692131552690 & 0.6411692131552690& 0 &0&0 &0\\
1.2537322752425418 &0.3905895060040396 &0.8631427692385082 & 0&0 &0\\
  1  &0.4274711580740817 &0.3555517808854274 & 0.21697706104049089&0 &0\\
  1  &0.3099153072147496 &0.3259623915325679 &-0.2881752086128284&0.6522975098655108 &0\\
  \hline
  &0.3099153072147496 &0.3259623915325679 &-0.2881752086128284&0.6522975098655108 &0\\
\end{tabular}

\vspace{0.15cm}
\begin{tabular}{c@{\hspace{0.16cm}}|c@{\hspace{0.16cm}}c@{\hspace{0.16cm}}c@{\hspace{0.16cm}}c@{\hspace{0.16cm}}c@{\hspace{0.16cm}}}
  0& 0& 0 & 0&0&0\\
0.641169213155269 &0.3031200089371227 &0.3380492042181466  &0 &0 &0\\
1.253732275242547&0.3905895060040396& 0.4629099915955034& 0.4002327776430044 &0 \\
1&0.4341539203752613&0.3418741772176282&0.2239719024071105 &0&0\\
  1&0.3099153072147496&0.3259623915325679&-0.2881752086128284&0&0.6522975098655108\\
  \hline
    &0.3099153072147496&0.3259623915325679&-0.2881752086128284&0&0.6522975098655108
\end{tabular}
\caption{Tableau for the explicit(top) and implicit(bottom) third order  $L$-stable  scheme.}
\label{tab:third order five stage}
  \end{table}
  The stability function is
  \begin{equation}
\begin{aligned}
&R (z)=\frac{-3.10127 z ^2-4.42559 z +11.3308}{- z ^3+6.98974 z ^2-15.7564 z +11.3308}.
\end{aligned}
\end{equation}
The method is L-stable.

Table \ref{tab:third order five stage 2} presents another five-stage, third-order method, stratifying $u_{n+1}=K_{s}$. It requires solving  four linear systems at each iteration, but more stable compare to the previous method . 
\begin{table}[h]
  \tiny
\centering
\begin{tabular}{c@{\hspace{0.16cm}}|c@{\hspace{0.16cm}}c@{\hspace{0.16cm}}c@{\hspace{0.16cm}}c@{\hspace{0.16cm}}c@{\hspace{0.16cm}}}
  0 & 0& 0 & 0&0 &0\\
0.3772977846271119 & 0.3772977846271119& 0 &0&0 &0\\
1 &0.3210924473454751 &0.6789075526545275 & 0&0 &0\\
  1  &0.2958359189953578 &0.3278679213986500 & 0.3762961596059923&0 &0\\
  1  &0.05826227065874467 &0.7093884017687849 &-0.2070619980550040&0.4394113256274744 &0\\
  \hline
  &0.05826227065874467 &0.7093884017687849 &-0.2070619980550040&0.4394113256274744 &0\\
\end{tabular}

\vspace{0.15cm}
\begin{tabular}{c@{\hspace{0.16cm}}|c@{\hspace{0.16cm}}c@{\hspace{0.16cm}}c@{\hspace{0.16cm}}c@{\hspace{0.16cm}}c@{\hspace{0.16cm}}}
  0& 0& 0 & 0&0&0\\
0.3772977846271117 &0.2709023139105694 &0.1063954707165423  &0 &0 &0\\
1&0.3210924473454735& 0.4580508073137827& 0.2208567453407465 &0 \\
1&0.4458748098646118&0.08691986121002987&0.3372847407465245 &0.1299205881788340&0\\
  1&0.05826227065874504&0.7093884017687844&-0.2070619980550035&-0.2178085843289785&0.6572199099564526\\
  \hline
    &0.05826227065874504&0.7093884017687844&-0.2070619980550035&-0.2178085843289785&0.6572199099564526
\end{tabular}
\caption{Tableau for the explicit(top) and implicit(bottom) third order  $L$-stable  scheme.}
\label{tab:third order five stage 2}
  \end{table}
  The stability function is
  \begin{equation}
\begin{aligned}
&R (z)=\frac{-35.1326 z^3-123.561 z^2-57.0133 z+498.399}{ z^4-23.1453 z^3+182.652 z^2-555.413 z+498.399}.
\end{aligned}
\end{equation}
The method is L-stable.

\section{Numerical Test}
\label{sec:numerical test}
In this section, we conduct a series of numerical tests for a scalar equation, a nonlinear diffusion equation, the Navier-Stokes equation and the Cahn-Hilliard equation, aiming to validate the convergence rate and  efficiency of the proposed method.

We now introduce the formulas and notation for the error metric. If the exact solution is available, we define the relative error of the algorithm as
\begin{equation}\label{eq:DefRelativeError}
E(h) := \frac{ \lVert I(h) - I_{E} \rVert_{\infty} }{ \lVert I_{E} \rVert_{\infty} },
\end{equation}
where $I(h)$ is the output of the method with time step size $h$, $\lVert \cdot \rVert_{\infty}$ denotes the infinity norm, and $I_{E}$ is the exact solution. If the exact solution is not available, we compute a reference solution using a high-order method with a fine time step size and treat it as the exact solution for error estimation.

With this notion, we define the numerical rate of convergence as
\begin{equation}\label{eq:DefConvergenceRate}
\mathcal{O}(h)=\log_r\frac{ E(rh) }{ E(h) }.
\end{equation}
where $r$ is the refinement factor. We set  it to be 2 in this section. 

Beyond analyzing convergence rates, we investigate the time step restrictions imposed by stability  if a steady solution exists and is an attractor for the initial conditions under consideration. Specifically, we numerically examine whether the solution of the initial value problem converges to a limiting function as time approaches infinity. For each method, we perform simulations with varying time step sizes and identify the largest step size for which the numerical solution remains stable and correctly approaches the limiting function. This extends the stability analysis from linear equations to nonlinear problems. In the implementation, we consider the numerical solution converged if the relative difference between the numerical solution and  the long-time limit solution is less than 0.01.

\subsection{Scalar equation}
We first test the proposed methods on a scalar ODE with a known exact solution to verify their convergence rate and assess whether they are well-conditioned:
\begin{equation}\label{eq:scalar equation}
\begin{aligned}
&y'=f (t,y)+G(t,y)y=\cos (t)y+ (- y+\cos(t))y.
\end{aligned}
\end{equation}
The exact solution is
\begin{equation}
\begin{aligned}
&y=\frac{e^{2 \sin (t)}}{1+ \int\limits_0^{t} e^{2 \sin (s)}\mathrm{d}s}.
\end{aligned}
\end{equation}
The results for the relative error and convergence rate are presented in Table~\ref{tab:scalar equation}. All methods resolve the solution to the machine precision level. The expected convergence rates are observed for all methods before the relative error reaches the machine precision.
\begin{table}[h]
\centering
\begin{tabular}{c|ccccccc}
\hline\hline
Method  &$E\left( \frac{1}{131072}\right)$ &  $\mathcal{O}$   &$E\left( \frac{1}{262144}\right)$ &  $\mathcal{O}$   &$E\left( \frac{1}{524288}\right)$ &  $\mathcal{O}$   &$E\left( 2^{-20}\right)$  \\ 
\hline\hline
Table \ref{tab:midpoint}&6.77e-13&2.00&1.69e-13&2.89&2.28e-14&0.80&1.31e-14\\ 
Table \ref{tab:trapezoid}&8.90e-13&2.14&2.03e-13&3.08&2.39e-14&1.29&9.75e-15\\ 
Table \ref{tab:second order scheme L}&1.73e-12&1.96&4.44e-13&1.53&1.53e-13&1.33&6.12e-14\\ 
Table \ref{tab:second order scheme 3}  &1.82e-12&1.97&4.67e-13&1.80&1.34e-13&2.42&2.50e-14\\
  \hline
  \hline
Method&$E\left( \frac{1}{1024}\right)$ &  $\mathcal{O}$   &$E\left( \frac{1}{2048}\right)$ &  $\mathcal{O}$   &$E\left( \frac{1}{4096}\right)$ &  $\mathcal{O}$   &$E\left( \frac{1}{8192}\right)$ \\  
\hline\hline
Table \ref{tab:third order four stage}&1.22e-12&3.03&1.49e-13&3.08&1.76e-14&2.35&3.46e-15\\
  Table \ref{tab:third order five stage}&1.42e-12&2.98&1.81e-13&3.52&1.57e-14&2.32&3.15e-15\\
   Table \ref{tab:third order five stage 2} &1.42e-12&2.98&1.81e-13&3.52&1.57e-14&2.32&3.15e-15\\ 
\hline
\end{tabular}
\caption{The relative errors and convergence rates of the solution to the scalar equation \eqref{eq:scalar equation} at $t=0.5$. The first column lists the corresponding table for each method. The first four methods are second order. The last third methods are third order.  }
\label{tab:scalar equation}
\end{table}

\subsection{Nonlinear diffusion equation}
A one-dimensional version of equation \eqref{eq:nonlinear advection-diffusion equation} is considered here,
\begin{equation}\label{eq:nonlinear advection-diffusion equation 1D}
\partial_{t} c = \partial_{x}\left( (1+\kappa c^{2}) \partial_{x} c \right) + S (t,x), \quad \left. c \right|_{t=0} = 0,
\end{equation}
on the domain $-\pi \leq x \leq \pi$ with periodic boundary conditions over the time interval $[0,1]$. The parameter $k$ is introduced to control the nonlinearity in this equation. Comparing to equation \eqref{eq:vector evolution equation}, we set $G(t, c) = \partial_{x}\left( (1+\kappa c^{2}) \partial_{x} \cdot \right)$ and $f (t,x,c) = S (t,x)$. 

The first goal of the numerical test is to examine the convergence rate of the proposed methods. To demonstrate their performance for a time-dependent source term, we set $S (t,x)=\cos(x)\sin(t)$ and $\kappa=1$. We approximate the spatial derivatives using a five-point finite difference stencil, constructed following the algorithm in \cite{fornberg1998classroom}. The equation is spatially discretized using 129 uniformly distributed grid points, and we vary the time step size to study convergence to the discretized system. Periodic boundary conditions are enforced by replacing the first and last rows of the matrix  $L=1 - a_{i,i} G (t_{n} + c_{i}h, K_{i-1}) K_{i}$ with vectors representing $c (\pi)-c (-\pi)=0$ and $c'(\pi)-c'(-\pi)=0$. 

Since the exact solution is not available, we compute a reference solution using the third-order method from Table~\ref{tab:third order four stage} with a time step size of $2^{-9} = \frac{1}{512}$. Given that the boundary conditions must be imposed explicitly, we employ the methods presented in Tables~\ref{tab:forward backward Euler}, \ref{tab:trapezoid}, \ref{tab:second order scheme L}, \ref{tab:third order five stage}, and \ref{tab:third order five stage 2}.  The results are summarized in Table~\ref{tab:nonlinear diffusion equation}, where the expected convergence rates are observed for all methods. To validate the reference solution, we recompute the table using an even smaller time step size of $2^{-10}$, and the results remain unchanged.
\begin{table}[h]
\centering
\begin{tabular}{c|ccccccc}
\hline\hline
Method  &$E\left( \frac{1}{16}\right)$ &  $\mathcal{O}$   &$E\left( \frac{1}{32}\right)$ &  $\mathcal{O}$   &$E\left( \frac{1}{64}\right)$ &  $\mathcal{O}$   &$E\left( \frac{1}{128}\right)$ \\ 
  \hline\hline
 Table \ref{tab:forward backward Euler}       &6.64e-02&1.00&3.33e-02&1.00&1.67e-02&1.00&8.33e-03\\ 
  \hline
Table \ref{tab:trapezoid} &9.49e-05&2.00&2.37e-05&2.00&5.91e-06&2.00&1.48e-06\\ 
Table \ref{tab:second order scheme L}&1.46e-04&1.98&3.70e-05&1.99&9.30e-06&2.00&2.33e-06\\ 
    \hline
  Table \ref{tab:third order five stage}    &1.35e-05&3.08&1.59e-06&3.00&1.99e-07&3.00&2.49e-08\\ 
 Table \ref{tab:third order five stage 2}&9.29e-06&2.88&1.26e-06&2.93&1.65e-07&2.98&2.09e-08\\ 
\hline
\end{tabular}
\caption{The relative errors and convergence rates of the solution to the nonlinear diffusion equation \eqref{eq:nonlinear advection-diffusion equation 1D}. The first column lists the corresponding table for each method. }
\label{tab:nonlinear diffusion equation}
\end{table}

Next, we aim to study the  time step restrictions imposed by stability. We set $S (t,x)=\cos x$. Then   equation \eqref{eq:nonlinear advection-diffusion equation 1D} admit the following long-time limit
\begin{equation}\label{eq:nonlinear diffusion equation long time limit}
\begin{aligned}
&c (\infty, x)=\frac{\sqrt[3]{2} \left(\sqrt{9 \kappa \cos ^2(x)+4}+3 \sqrt{\kappa} \cos (x)\right)^{2/3}-2}{2^{2/3} \sqrt{\kappa} \sqrt[3]{\sqrt{9 \kappa \cos ^2(x)+4}+3 \sqrt{\kappa} \cos (x)}}.
\end{aligned}
\end{equation}
For an explicit RK method, since there is no additional mechanism to enforce the boundary condition, the numerical solution exists only for a short time before eventually blowing up, regardless of how small the time step size is.
Therefore, we consider the IMEX method instead. One way to apply an IMEX method while solving only a linear system is to treat the linear part of the operator implicitly, a strategy we refer to as IMEX with linear splitting. In this case, we treat $ \partial_{x}^{2} c$ implicitly while handling $ \partial_{x} \left(\kappa c^{2} \partial_{x} c \right) + S (t,x)$ explicitly. Although the presence of the second-order derivative in the explicit term still imposes a time step restriction of $\Delta t \sim \Delta x^{2}$, this approach allows larger time steps compared to treating the entire diffusion term explicitly, especially when $\kappa$ is small. Moreover, it ensures that the boundary condition is satisfied at each step. In the following numerical tests, we employ the second-order IMEX method described in Section 2.5 of \cite{ascher1997implicit}.

Table \ref{tab:nonlinear diffusion equation time step size} presents the estimated upper bound of the time step size required for the numerical solution to converge to the theoretical long-time limit. When $\kappa$ is large, the time step size for IMEX with linear splitting is approximately $10^{3}$ times smaller than that of the proposed semi-IMEX method. When $\kappa$ is small, the time step size for IMEX with linear splitting remains smaller than that of the semi-IMEX method, but their orders of magnitude are comparable. The forward-backward Euler scheme (shown in Table \ref{tab:forward backward Euler}) appears to be unconditionally stable in this case, as the solution converges even with a time step size of $10^{4}$. For all semi-IMEX methods, the time step size scales inversely with $\kappa$.

Among the second order method, the L-stable method (Table \ref{tab:second order scheme L}) performs better than the A-stable method (Table \ref{tab:trapezoid}). For the third-order methods, the one in Table \ref{tab:third order five stage} has a smaller stable time step size, possibly because it solves three linear systems per iteration, whereas the method in Table \ref{tab:third order five stage 2} solves four.

Finally, we note that the results also depend on the initial condition. The result in table \ref{tab:nonlinear diffusion equation time step size} is for the zero initial value. If the initial condition is closer to the long-time limit \eqref{eq:nonlinear diffusion equation long time limit}, all semi-IMEX methods permit a larger time step size, whereas the step size for IMEX remains relatively unchanged.
\begin{table}[h]
\centering
\begin{tabular}{c|cccccc}
\hline\hline
$\kappa$  & Table \ref{tab:forward backward Euler}    & Table \ref{tab:trapezoid}   &Table \ref{tab:second order scheme L} & Table \ref{tab:third order five stage}   &Table \ref{tab:third order five stage 2}  & IMEX \\ 
  \hline\hline
  0.25&$>10^{4}$&27.5&117&5.86&16.3&2.73\\
  0.5&$>10^{4}$&14.0&23.4&3.42&8.60&0.033\\
  1&$>10^{4}$&4.59&9.52&2.14&5.60&0.0068\\
  2&$>10^{4}$&2.13&3.91&1.29&3.23&0.0033\\
  4&$>10^{4}$&1.14&1.93&0.891&1.95&0.0019\\
\hline
\end{tabular}
\caption{ The estimated upper bound of the time step size required for the numerical solution to converge to the theoretical long-time limit of the nonlinear diffusion equation \eqref{eq:nonlinear advection-diffusion equation 1D}. The first column lists the values of $\kappa$, and the last column (IMEX) refers to IMEX with linear splitting. }
\label{tab:nonlinear diffusion equation time step size}
\end{table}

\subsection{Cahn-Hilliard equation}
The Cahn-Hilliard equation models phase separation in binary mixtures, describing the evolution of two components over time to form distinct regions, which takes form
\begin{equation}\label{eq:Cahn-Hilliard equation}
\partial_{t} \phi = \nabla \cdot \left( M \nabla \mu \right),\quad  \mu = -\epsilon^2 \Delta \phi + f'(\phi),
\end{equation}
where \( \phi \) is the phase field, \( \mu \) is the chemical potential, \( M \) is the mobility, and \( f(\phi)=\frac{1}{4}(\phi^2 - 1)^2 \) is the double-well potential. Here, we set $M=1$ for simplicity. For a bounded domain, the no-flux boundary conditions are typically applied \cite{miranville2019cahn}:
\begin{equation}
\begin{aligned}
& \partial_{\mathbf{n}}\phi=0, \quad \partial_{\mathbf{n}}\mu=\partial_{\mathbf{n}}\left(  -\epsilon^2 \Delta \phi + f'(\phi) \right)= 0,
\end{aligned}
\end{equation}
where $\mathbf{n}$ is the normal vector of the boundary. To match the form of equation \eqref{eq:vector evolution equation}, we define $G (\phi)= -\epsilon^{2}\Delta^{2}\cdot + \nabla\cdot \left(    \left( 3 \phi ^2-1 \right) \nabla \cdot \right)$.  Although  the boundary condition is nonlinear in $\phi$, the proposed method still applies. We can rewrite the boundary condition as $B (\phi)\phi=0$, where the operator is defined as $B (\phi)= (\partial_{\mathbf{n}},  -\epsilon^{2}\partial_{\mathbf{n}} \Delta  + f'' (\phi)\partial_{\mathbf{n}} )$. In the implementation of the semi-IMEX RK method,  we first construct the matrix  $L=1 - a_{i,i} G (t_{n} + c_{i}h, K_{i-1}) K_{i}$ and then replace some rows of matrix $L$ with $B (\phi)$. 

For simplicity, we consider the one dimensional case with the  initial condition  $\left. \phi \right|_{x=0}=\tanh (x)$,  $\epsilon=1$. The domain is $[-20, 20]$ and the time interval is $[0,1]$. A non-uniform mesh with $128$ points is used, with grid points denser around $x=0$  and sparser near the endpoints.

First, we aim to examine the convergence rate of the proposed methods. The results are shown in Table~\ref{tab:Cahn-Hilliard equation}, where the expected convergence rates are observed for all methods, showing the validity of proposed method for the problem with nonlinear boundary conditions. Since the exact solution is not available, we compute a reference solution using the third-order method from Table~\ref{tab:third order five stage 2} with time step sizes of $2^{-13} = \frac{1}{8192}$ and $2^{-14}$. The results in Table~\ref{tab:Cahn-Hilliard equation} are the same for both reference solutions, confirming convergence.
\begin{table}[h]
\centering
\begin{tabular}{c|ccccccc}
\hline\hline
Method&$E\left( \frac{1}{256}\right)$ &  $\mathcal{O}$   &$E\left( \frac{1}{512}\right)$ &  $\mathcal{O}$   &$E\left( \frac{1}{1024}\right)$ &  $\mathcal{O}$   &$E\left( \frac{1}{2048}\right)$ \\ 
  \hline\hline
  Table \ref{tab:forward backward Euler}  &8.41e-05&1.00&4.20e-05&1.00&2.10e-05&1.00&1.05e-05\\ 
  \hline
Table \ref{tab:trapezoid}&2.32e-07&1.94&6.05e-08&1.96&1.55e-08&1.98&3.94e-09\\ 
Table \ref{tab:second order scheme L}&2.04e-07&1.99&5.14e-08&1.99&1.29e-08&2.00&3.23e-09\\ 
  \hline
  Table \ref{tab:third order five stage}  &5.52e-08&2.78&8.03e-09&2.86&1.10e-09&2.94&1.43e-10\\ 
  Table \ref{tab:third order five stage 2} &3.07e-08&2.97&3.91e-09&3.07&4.68e-10&3.24&4.95e-11\\ 
\hline
\end{tabular}
\caption{Relative errors and convergence rates of the solution to the Cahn-Hilliard equation \eqref{eq:Cahn-Hilliard equation}.}
\label{tab:Cahn-Hilliard equation}
\end{table}

Second, we study the time step restrictions imposed by stability. For the unbounded domain, the steady solution of equation \eqref{eq:Cahn-Hilliard equation} is  $c=\tanh \left( \frac{x}{\sqrt{2}\epsilon} \right)$ for the initial condition  $\left. \phi \right|_{x=0}=\tanh (x)$. Since the domain is bounded, we numerically solve the steady solution using Newton's method. The IMEX method with linear splitting treats $-\epsilon^{2} \Delta^{2}\phi$ implicitly and $\Delta f' (\phi)$ explicit. The time step size restrictions for different methods are presented in Table \ref{tab:Cahn-Hilliard equation time step size}.

The stable time step size for the IMEX method with linear splitting is smaller than that for the semi-IMEX method, differing by a factor of 3 to 66 when $\epsilon = 1$. However, this difference is much smaller compared to the case of the nonlinear diffusion equation, where the factor is around 1000. This discrepancy arises because inverting $-\epsilon^{2} \Delta^{2}$ introduces a $k^{-4}$ weighting for high-wavenumber modes, which dampens the $k^{2}$ weight added by the explicit part. As $\epsilon$ decreases, this damping effect weakens, leading to a larger difference in stable time step sizes. When $\epsilon = 0.25$, the stable time step size for the IMEX method with linear splitting differs from that of the semi-IMEX method by a factor of 10 to 100, reflecting this reduced damping effect. Consistent with the results for the nonlinear diffusion equation, the methods listed in Tables \ref{tab:forward backward Euler}, \ref{tab:second order scheme L}, and \ref{tab:third order five stage 2} exhibit better performance in terms of time step size restrictions. 

\begin{table}[h]
\centering
\begin{tabular}{c|cccccc}
\hline\hline
$\epsilon$  & Table \ref{tab:forward backward Euler}    & Table \ref{tab:trapezoid}   &Table \ref{tab:second order scheme L} & Table \ref{tab:third order five stage}   &Table \ref{tab:third order five stage 2}  & IMEX \\ 
  \hline\hline
    0.25&1.17&0.198&0.452&0.111&0.625&0.011\\
  0.5&6.72&0.745&1.83&0.334&2.49&0.08\\
  1&24.6&2.12&7.72&1.10&8.45&0.36\\
\hline
\end{tabular}
\caption{ The estimated upper bound of the time step size required for the numerical solution to converge to the steady solution of the Cahn-Hilliard equation \eqref{eq:Cahn-Hilliard equation}. The first column lists the values of $\epsilon$, and the last column (IMEX) refers to IMEX with linear splitting. }
\label{tab:Cahn-Hilliard equation time step size}
\end{table}

\subsection{Navier-Stokes equation}
The pressure-driven flow of an incompressible, homogeneous fluid is governed by the Navier-Stokes equations:
\begin{equation}\label{eq:NSE}
\begin{aligned}
& \partial_{t} \mathbf{u} + \mathrm{Re} \mathbf{u} \cdot \nabla \mathbf{u} =  \Delta \mathbf{u} -\nabla p + f \mathbf{e}_{x}, \quad \nabla \cdot \mathbf{u} = 0,
\end{aligned}
\end{equation}
where $\mathbf{u} = (u,v)$ is the velocity field, $p$ is the pressure field and $\mathrm{Re}$ is the Reynolds number. $f$ is a constant and the term $f \mathbf{e}_{x}$ represents the imposed pressure gradient in the $x$-direction.

We consider this flow in a channel with non-flat boundaries, defined as
\begin{equation}
\left\{ (x,y) \mid -1<x<1, h_{1} (x) \leq y \leq h_{2}(x) \right\}.
\end{equation} No-slip boundary conditions are enforced at the walls $y = h_1(x)$ and $y = h_2(x)$, while periodic boundary conditions are applied at the inlet and outlet $x=\pm1$. Such a domain has numerous applications in microfluidic devices \cite{watanabe2007simplified,durst1993plane}. By appropriately choosing $h_{1}$ and $h_{2}$, the channel can be tailored for specific applications, such as microcentrifuges \cite{khojah2017size,haddadi2017inertial,che2016classification}. These devices selectively trap large particles in flowing suspensions, enabling applications like isolating  tumor cells from smaller red blood cells in patient blood samples.

The IMEX method typically treats the advection term explicitly and the viscosity term implicitly. For example, it can be written as
\begin{equation}\label{eq:NS Stokes}
\begin{aligned}
& \frac{\mathbf{u}^{n+1}- \mathbf{u}^{n}}{h} +\mathrm{Re}  \mathbf{u}^{n} \cdot \nabla \mathbf{u}^{n}  = \Delta \mathbf{u}^{n+1} -\nabla p^{n+1}  + f \mathbf{e}_{x}, \quad \nabla \cdot \mathbf{u}^{n+1} = 0.
\end{aligned}
\end{equation}
It is equivalent to solve a generalized Stokes equation for $\mathbf{u}^{n+1}$ with a forcing term depending on $\mathbf{u}^{n}$ and boundary conditions. This is a  standard method discussed in many tutorial such as  \cite{hecht2005freefem++}.  We propose an alternative formulation by considering the advection term partially implicitly. Under this framework, the first-order semi-IMEX method \eqref{eq:forward-backward Euler} leads to the following scheme:
\begin{equation}\label{eq:NS Oseen}
\begin{aligned}
& \frac{\mathbf{u}^{n+1}- \mathbf{u}^{n}}{h} +\mathrm{Re} \mathbf{u}^{n} \cdot \nabla \mathbf{u}^{n+1}  =  \Delta \mathbf{u}^{n+1} -\nabla p^{n+1}+ f \mathbf{e}_{x}, \quad \nabla \cdot \mathbf{u}^{n+1} = 0.
\end{aligned}
\end{equation}
It is a generalized Oseen equation for $\mathbf{u}^{n+1}$ with a forcing term that depends on $\mathbf{u}^{n}$. It is well known that the Stokes equation is not a uniform asymptotic approximation of the Navier-Stokes equation for small Reynolds numbers, whereas the Oseen equation is \cite{happel2012low}. Furthermore, as $h \to \infty$, scheme \eqref{eq:NS Oseen} reduces to the Picard iteration, which is commonly used for solving the steady Navier-Stokes equation due to its stability and global convergence properties  \cite{pollock2019anderson}.  In general, scheme \eqref{eq:NS Oseen} allows for a much larger time step size compared to scheme \eqref{eq:NS Stokes}.

The method in table \ref{tab:trapezoid} gives us a second order accuracy in time scheme
\begin{equation}\label{eq:NSE trapezoid}
\begin{aligned}
 & \frac{K_2- \mathbf{u}^{n}}{\frac{1}{2}h} + \mathrm{Re}\mathbf{u}^{n} \cdot \nabla K_2  =  \Delta K_2 -\nabla p_{2}^{n+1}+ f \mathbf{e}_{x}, \quad \nabla \cdot K_2 = 0,\\
 & \frac{K_3- \mathbf{u}^{n}}{\frac{1}{2}h} +\mathrm{Re}  K_{2} \cdot \nabla K_3  = \Delta K_3 -\nabla p_{3}^{n+1} + f \mathbf{e}_{x}, \quad \nabla \cdot K_3 = 0,\\
 &\mathbf{u}^{(n+1)}=2K_{3}-\mathbf{u}^{n}. 
\end{aligned}
\end{equation}
We prefer this method for its simplicity of implementation. Given a solver for equation \eqref{eq:NS Oseen}, achieving a second-order accurate scheme in time requires only two applications of the solver with different velocity fields. Moreover, since the no-slip and periodic boundary conditions and incompressibility are enforced when solving for $K_{2}$ and $K_{3}$, we ensure that $\mathbf{u}^{(n+1)}$ satisfies the boundary condition and incompressibility.

We use the finite element method to spatially discretize the Navier-Stokes equation. The algorithm is implemented using the open-source package FreeFEM++ \cite{hecht2012new}. A triangular mesh with second-order polynomial basis functions (P2) is employed. The maximum edge length of the triangular mesh is set to 0.06. We set the initial condition for the velocity field to zero. For the domain, we define $h_{1}(0)=0$ and $h_{2}(x) = \frac{1}{2} \tanh \left( 10 \cos \left( \frac{4 \pi x}{5} \right) \right) + \frac{3}{2}$. The time interval is $[0,1]$. Since the exact solution is unavailable, the reference solution is computed using a time step size of $h = \frac{1}{128}$ and the scheme \eqref{eq:NSE trapezoid} . The relative error and convergence rate of the scheme in \eqref{eq:NSE trapezoid} are presented in Table \ref{tab:NSE trapezoid} for the parameters $\mathrm{Re}=100, f=-1$. The second order convergence rate for the scheme in \eqref{eq:NSE trapezoid} is obtained. 
\begin{table}[h]
\centering
\begin{tabular}{c|ccccccccc}
\hline\hline
Method   &$E\left( \frac{1}{2}\right)$ &  $\mathcal{O}$   &$E\left( \frac{1}{4}\right)$ &  $\mathcal{O}$   &$E\left( \frac{1}{8}\right)$ &  $\mathcal{O}$  &$E\left( \frac{1}{16}\right)$ &  $\mathcal{O}$   &$E\left( \frac{1}{32}\right)$  \\ 
  \hline\hline
Eq \eqref{eq:NSE trapezoid} &3.70e-01&2.17&8.20e-02&1.96&2.11e-02&2.02&5.21e-03&2.10&1.21e-03\\ 
\hline
\end{tabular}
\caption{Relative errors and convergence rates of the solution to the Navier Stokes equation \eqref{eq:NSE}.}
\label{tab:NSE trapezoid}
\end{table}

Next, we examine the time step restrictions imposed by stability. The steady solution to the Navier-Stokes equation \eqref{eq:NSE} is computed using Newton's method. The schemes in Eq. \eqref{eq:NS Oseen} and Eq. \eqref{eq:NSE trapezoid} appear to be unconditionally stable in this case, as the solution converges even with a time step size of $10^{4}$. In contrast, the first-order IMEX method in Eq. \eqref{eq:NS Stokes} requires a much smaller time step size, and the required step size decreases as $\mathrm{Re}$ increases, indicating increasingly stringent stability constraints for higher Reynolds numbers. The absence of time step restrictions for Eq. \eqref{eq:NS Oseen} and Eq. \eqref{eq:NSE trapezoid} suggests that both methods are well suited for high-$\mathrm{Re}$ flows where large time steps are desirable.
\begin{table}[h]
\centering
\begin{tabular}{c|ccc}
\hline\hline
$\mathrm{Re}$  &Eq \eqref{eq:NS Oseen}&  Eq \eqref{eq:NSE trapezoid} &Eq \eqref{eq:NS Stokes}  \\ 
  \hline\hline
  100&$>10^{4}$ &$>10^{4}$  & 0.048\\
  200&$>10^{4}$  &$>10^{4}$  & 0.0048 \\
  300&$>10^{4}$  &$>10^{4}$  & 0.0018 \\
\hline
\end{tabular}
\caption{The estimated upper bound of the time step size required for the numerical solution to converge to the Navier Stokes equation \eqref{eq:NSE}.  }
\label{tab:NSE time step size}
\end{table}

\subsection{Other applications}
Another important application of this semi-IMEX RK scheme is in inhomogeneous fluids, where properties such as density, viscosity, or composition vary spatially. In many cases, this variation depends on the concentration field of a solute advected by the fluid, such as in ocean and atmosphere \cite{more2023motion, ecke2023turbulent, ding2023dispersion}. The governing equations are given by:
\begin{equation}
\begin{aligned}
& \rho(c) \left( \partial_t \mathbf{u} + \mathbf{u} \cdot \nabla \mathbf{u} \right) = \mu(c) \Delta \mathbf{u} - \nabla p + \rho(c) g, \quad \nabla \cdot \mathbf{u} = 0, \
& \partial_t c + \mathbf{u} \cdot \nabla c = \kappa \Delta c.
\end{aligned}
\end{equation}
Here, $c$ is the concentration field, $\kappa$ is the diffusivity of the solute, $\rho(c)$ is the fluid density, and $\mu(c)$ is the dynamic viscosity, both of which depend on $c$. When the variation in density or viscosity is small, IMEX with linear splitting is sufficient. However, when the variation is large, the IMEX method with linear splitting requires very small time step sizes due to stability constraints, whereas the proposed semi-IMEX RK scheme helps alleviate this restriction. Since this system shares characteristics with both the nonlinear diffusion equation and the Navier-Stokes equations, a detailed numerical investigation is omitted here.

 \section{Conclusion}
 \label{sec:conclusion}
We have proposed a novel semi-IMEX Runge-Kutta method for solving ODEs with both stiff and non-stiff components, particularly when the stiff term can be expressed as a matrix-vector product. Such ODEs arise from PDEs with nonlinear diffusion terms. By treating the matrix in the stiff term explicitly while keeping the solution variable implicit, our approach significantly reduces computational complexity by eliminating the need to solve nonlinear systems or compute Jacobians.

We have developed a class of semi-IMEX RK methods with varying orders of accuracy and analyzed their stability and convergence properties. Numerical experiments on scalar equations, nonlinear diffusion models, the Navier-Stokes equations, and the Cahn-Hilliard equation demonstrate the efficiency and accuracy of the proposed scheme. Compared to the IMEX-RK method with linear splitting, our approach allows for larger time step sizes while still satisfying stability requirements, resulting in substantial improvements in computational cost.

The repository containing the implementation of the proposed Runge-Kutta methods and scripts to reproduce selected tables is available at \verb|https://github.com/wdachub/semi-IMEX-RK-method|.

Future research directions include several avenues. First, the proposed modification to equation \eqref{eq:IMEX-RK scheme} in equation \eqref{eq:semi-IMEX-RK scheme} introduces additional coupling when deriving the order condition for semi-IMEX methods, making the development of higher-order methods, such as a fourth-order scheme, challenging. A clever method for deriving the order conditions or exploring alternative modifications to the classical IMEX method will be explored in future work. Second, for nonlinear equations, other stability requirements, such as $B$-stability \cite{butcher1975stability}, will also be studied in future research.

\section{ Declaration of Interests}
The author report no conflict of interest.

\end{document}